\documentclass[11pt,a4paper]{article}
\usepackage[left=2cm, top=2.5cm,bottom=2.5cm,right=2cm]{geometry}
\usepackage{mathtools,amssymb,amsthm,mathrsfs,calc,graphicx,xcolor,cleveref,dsfont,tikz,pgfplots,bm,url,tabularx}
\usepackage[british]{babel}
\usepackage{amsfonts}              
\usepackage[T1]{fontenc}
\usepackage{enumitem}
\usepackage[color=green]{todonotes}
\usepackage{tikz-cd}
		\usetikzlibrary{calc} 
		\usetikzlibrary{intersections}
\usepackage[labelfont=bf]{caption}
\usepackage[font=small]{subcaption}

\theoremstyle{plain}
\newtheorem{theorem}{Theorem}

\newtheorem{corollary}[theorem]{Corollary}

\theoremstyle{definition}

\theoremstyle{remark}

\newcommand{\dint}{\textup{d}}

\def\EE{\mathbb{E}}

\def\PP{\mathbb{P}}

\def\RR{\mathbb{R}}

\makeatletter
\let\@fnsymbol\@alph
\makeatother

\begin{document}

\title{\bfseries The proportion of triangles in a class of\\ anisotropic Poisson line tessellations}

\author{Nils Heerten\footnotemark[1],\; Julia Krecklenberg\footnotemark[2],\; Christoph Th\"ale\footnotemark[3]}

\date{}
\renewcommand{\thefootnote}{\fnsymbol{footnote}}
\footnotetext[1]{Ruhr University Bochum, Germany. Email: nils.heerten@rub.de}
\footnotetext[2]{Ruhr University Bochum, Germany. Email: julia.krecklenberg@rub.de}
\footnotetext[3]{Ruhr University Bochum, Germany. Email: christoph.thaele@rub.de}

\maketitle

\begin{abstract}
\noindent
Stationary Poisson processes of lines in the plane are studied whose directional distributions are concentrated on $k\geq 3$ equally spread directions. The random lines of such processes decompose the plane into a collection of random polygons, which form a so-called Poisson line tessellation. The focus of this paper is to determine the proportion of triangles in such tessellations, or equivalently, the probability that the typical cell is a triangle. As a by-product, a new deviation of Miles' classical result for the isotropic case is obtained by an approximation argument. \\

\end{abstract}

\section{Introduction and results}

\begin{figure}[t]
	\begin{minipage}{0.33\textwidth}
		\includegraphics[width=0.9\columnwidth]{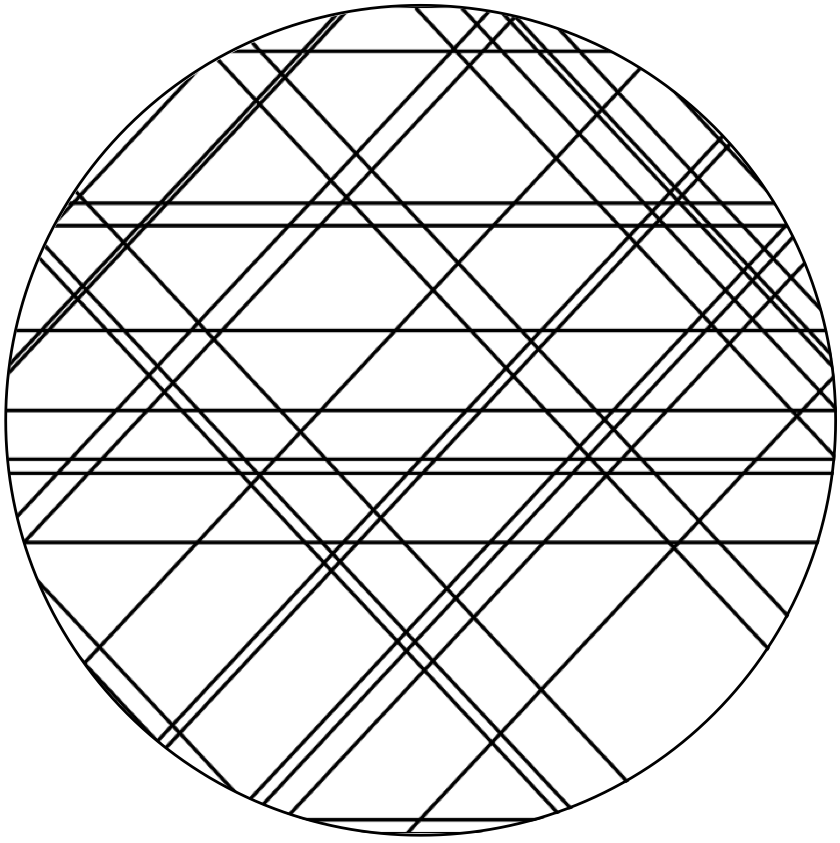}
	\end{minipage}%
	\begin{minipage}{0.33\textwidth}
		\centering
		\includegraphics[width=0.9\columnwidth]{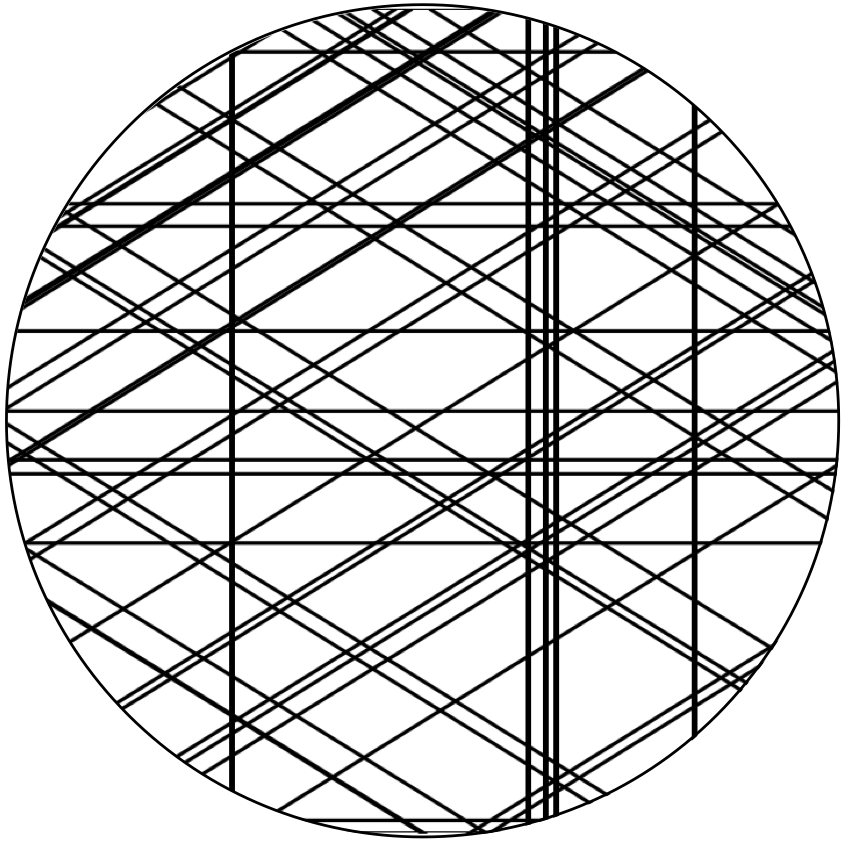}
	\end{minipage}
	\begin{minipage}{0.33\textwidth}
		\centering
		\includegraphics[width=0.9\columnwidth]{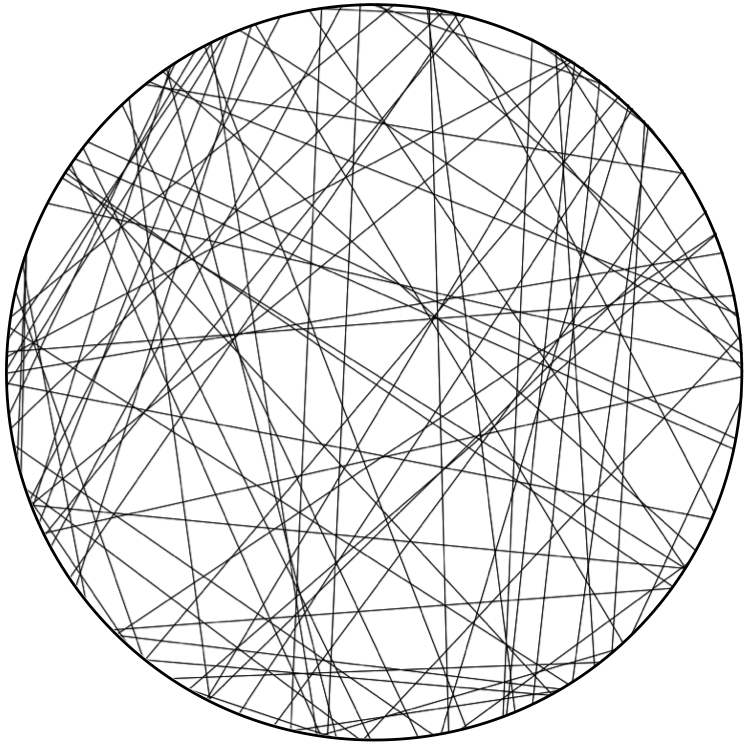}
	\end{minipage}%
	\caption{Simulation of a Poisson line tessellation with directional distribution $G_3$ (left), $G_4$ (middle) and $G_{\rm unif}$ (right).}
	\label{fig:PLP}
\end{figure}

The study of random polygons induced by a Poisson process of random lines in the plane is among the most classical topics in stochastic geometry. The distribution of a stationary Poisson line process $X=X(\gamma,G)$ in the plane is completely determined by its intensity $\gamma>0$ and its directional distribution $G$. For us, the latter is a probability measure on the interval $[0,\pi)$ satisfying $G(\theta)<1$ for each $\theta\in[0,\pi)$. We refer to the monograph of \cite{SKM} for further background material and detailed descriptions and explanations. The typical cell $Z=Z(\gamma,G)$ of a stationary Poisson line tessellation with intensity $\gamma$ and directional distribution $G$ can intuitively be thought of as a random polygon selected `uniformly at random' among the collection of all polygons (in a very large observation window) induced by $A$, regardless of size and shape. Formally, its distribution can be defined using Palm calculus as explained in detail in \cite{SKM}, see also \eqref{eq:TypicalCellDistributio} below.

In this paper we are interested in the probability that the typical cell $Z$ is a triangle. Since the intensity $\gamma$ only acts as a scaling parameter, this probability cannot depend on $\gamma$ and we can take $\gamma=1$ for simplicity and write $Z(G)$ instead of $Z(1,G)$. Further, we define the triangle probability
$$
p_3(G) := \PP[\,Z(G)\text{ is a triangle}\,],
$$
which can equivalently be described as the proportion of triangles among the polygons of the Poisson line tessellation:
$$
p_3(G) = \lim_{R\to\infty}{1\over \pi R^2}{\EE}\sum_{c\subset B_R}{\bf 1}\{c\text{ is a triangle}\},
$$
where $B_R$ stands for a circle of radius $R>0$ centred at the origin and the sum runs over all tessellation cells $c$ contained in $B_R$.
If the directional distribution $G=G_{\rm unif}$ is the uniform distribution on $[0,\pi)$ and the Poisson line tessellation is isotropic, it is known since \cite{Miles64} (see Theorem 6 therein) that
\begin{equation}\label{eq:p3Uniform}
	p_3(G_{\rm unif}) = 2-{\pi^2\over 6}\approx 0.35507,
\end{equation}
compare also with \cite{Miles73} and with the computations given in Section \ref{sec:IsotropicCase}. A  realization of an isotropic Poisson line process is shown in the right panel of Figure~\ref{fig:PLP}. In Section~\ref{sec:alt_miles} we will provide an alternative proof for \eqref{eq:p3Uniform} using new results from the present paper. We further remark that in the isotropic case also the probability
\begin{align*}
	&\PP[\,Z(G_{\rm unif})\text{ is a quadrangle}\,]= \pi^2\log 2-{1\over 3}-{7\pi^2\over 36}-{7\over 2}\sum_{i=1}^\infty{1\over i^3}\approx 0.381466
\end{align*}
is known from \cite{Tanner83}. However,  for $k\geq 5$ the probabilities $\PP[\,Z(G_{\rm unif})\text{ has exactly $k$ vertices}\,]$ can be expressed only as rather involved multiple integrals, which can be evaluated numerically, see \cite{Calka03}. On the other hand, it is well known that the expected number of vertices of the typical cell is $4$, independently of the choice of the directional distribution $G$, see \cite[Section 10.5.1]{SKM}.

On the other extreme, if $G$ is concentrated on only two different values, all cells are almost surely parallelograms. So, in this case we have $p_3(G)=0$. Thus, the next non-trivial case arises if the directional distribution $G$ is concentrated on three different values. For simplicity and concreteness we focus on the case where $G$ is given by
\begin{align}\label{eq:G3pq}
	G_3(p,q):=p\delta_0+q\delta_{\pi\over 3}+(1-p-q)\delta_{2\pi\over 3},
\end{align}
where we write $\delta_{(\cdot)}$ for the Dirac measure and where $p,q\in(0,1)$ are weights satisfying $0<p+q<1$.
In other words, $G_3(p,q)$ is concentrated on the angles $0$, $\pi/3$ and $2\pi/3$ with weights $p$, $q$ and $1-p-q$, respectively. A simulation of a Poisson line tessellation with directional distribution $G_3(1/3,1/3)$ is shown in the left panel of Figure \ref{fig:PLP}. We remark that a stationary Poisson line process with directional distribution $G_3(1/3,1/3)$ is of course not invariant under \textit{all} rotations in the plane. However, it is invariant under rotations whose angle is an integer multiple of $\pi/3$. The corresponding Poisson line tessellation can thus be called $G_3(1/3,1/3)$-pseudo isotropic.

\begin{figure}[t]
	\begin{center}
		\includegraphics[width=0.5\columnwidth]{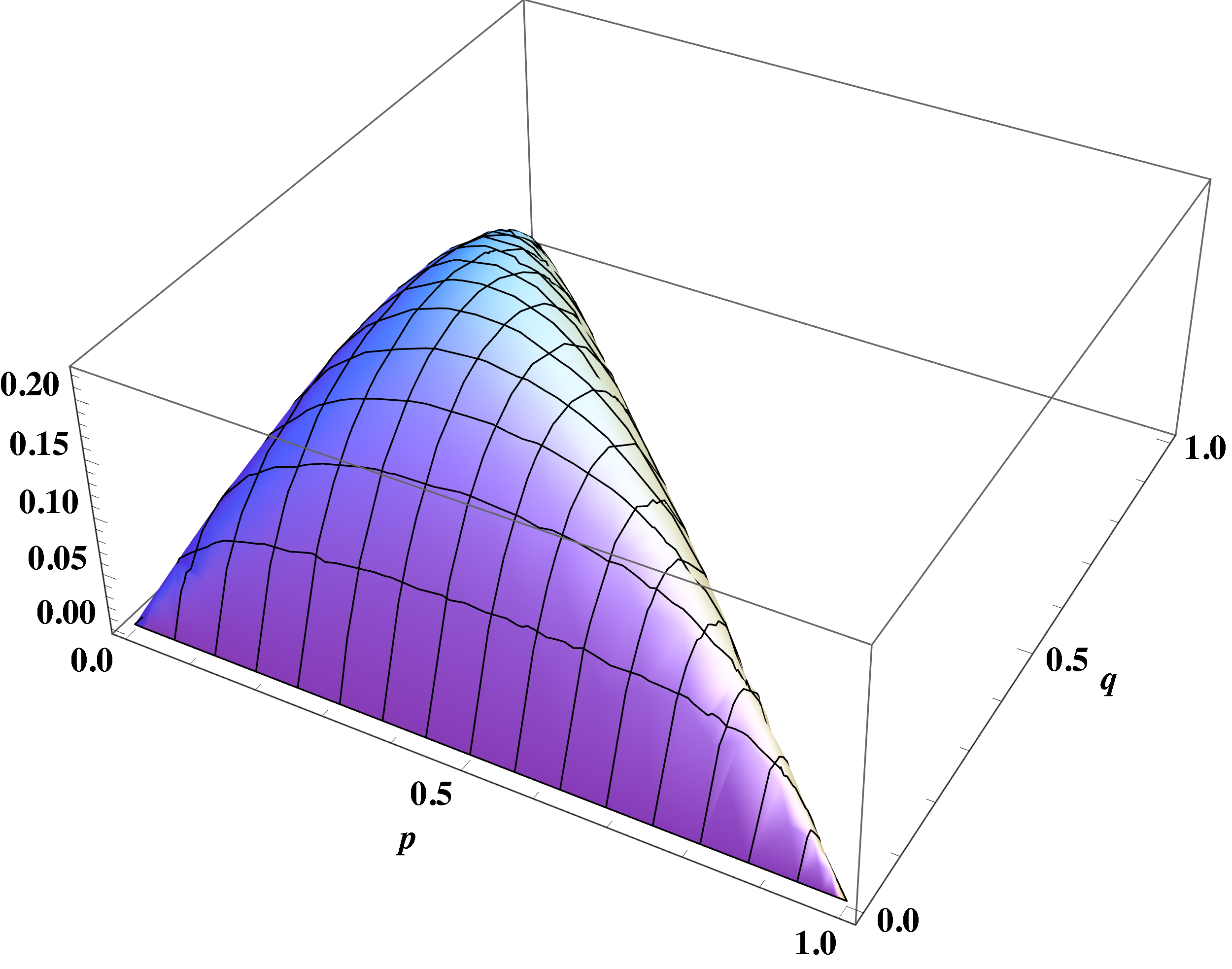}
	\end{center}
	\caption{A plot of $p_3(G_3(p,q))$}
	\label{fig:G3}
\end{figure}

Our first result is a formula for $p_3(G_3(p,q))$ in terms of the weights $p$ and $q$. Also, we determine those weights for which $p_3(G_3(p,q))$ attains its maximal value, see Figure \ref{fig:G3}.

\begin{theorem}\label{thm:3richtungen}
	For all $0<p,q<1$ with $0<p+q<1$, we have that
	$$
	p_3(G_3(p,q)) ={2pq(1-p-q)\over p+q-p^2-q^2-pq}. 
	$$
	The maximal value for $p_3(G_3(p,q))$ is attained precisely if $p=q=1/3$ and is given by
	$$
	\max_{0<p+q<1}p_3(G_3(p,q))  = p_3(G_3(1/3,1/3))  = {2\over 9}.
	$$
\end{theorem}

In analogy with the case of three directions just studied, one can consider a Poisson line tessellation with directional distribution $G_4(p,q,r):=p\delta_0+q\delta_{\pi\over 4}+r\delta_{\pi\over 2}+(1-p-q-r)\delta_{3\pi\over 4}$ with weights $0<p,q,r<1$ satisfying $0<p+q+r<1$, as shown in the middle panel of Figure~\ref{fig:PLP}. The corresponding triangle probability is in this case given by
\begin{align*}
	p_3(G_4(p,q,r))=&\frac{2 p}{\sqrt{2}p+2q+\sqrt{2}r-\sqrt{2}p^2-2q^2-\sqrt{2}r^2-2pq+\left(2-2\sqrt{2}\right)pr-2qr}\\
	&\times\bigg(\frac{3 q r}{{2+p(-2+\sqrt{2})-q+r(-2+\sqrt{2})}}+\frac{3\sqrt{2} q(1-p-q-r)}{2+(-2+\sqrt{2}) p+r(-2+2\sqrt{2})}\\
	&\phantom{{}=\cdot\Bigg({}}+\frac{2r\sqrt{2}(1-p-q-r)}{\sqrt{2}+p(2-\sqrt{2})+\sqrt{2}q+r(2-\sqrt{2})}\bigg),
\end{align*}
as demonstrated in \cite{Krecklenberg}. Since the triangle probabilities for five or more directions with arbitrary weights become increasingly more involved, from now on we concentrate on the special case where all weights are equal. Namely, we take for integers $k\geq 3$ as directional distribution the probability measure
$$
G_k := {1\over k}\sum_{\ell=0}^{k-1}\delta_{\ell\pi\over k},
$$
which for $k=3$ and $k=4$ reduces to $G_3(1/3,1/3)$ and $G_4(1/4,1/4,1/4)$, respectively. In other words, $G_k$ puts weight $1/k$ onto $k$ equally spread directions. The Poisson line tessellation induced by such a directional distribution is $G_k$-pseudo isotropic in that it is invariant under rotations in the plane whose angle is an integer multiple of $\pi/k$. In our second result we determine the triangle probabilities $p_3(G_k)$.

\begin{theorem}\label{thm:krichtungen}
	For $k\geq 3$ we have that
	\begin{align*}
		p_3(G_k) &= {4\over k}\tan^2{\pi\over 2k}\sum_{i=1}^{k-2}\Bigg[(k-i)\sum_{j=1}^{k-i-1}{\sin{i\pi\over k}\,\sin{j\pi\over k}\,\sin{(i+j)\pi\over k}\over\sin{i\pi\over k}+\sin{j\pi\over k}+\sin{(i+j)\pi\over k}}\Bigg].
	\end{align*}
\end{theorem}

\begin{figure}[t!]
	\begin{minipage}[b]{0.5\textwidth}
		\centering
		{\footnotesize
			\begin{tabular}{c|c|c|c}
				$k=3$ & $k=4$ & $k=5$ & $k=6$\\
				\\[-1em]
				\hline
				\\[-1em]
				$2\over 9$ & $4(5\sqrt{2}-7)$ & ${32\over\sqrt{5}}-14$ & ${4\over 3}(70\sqrt{3}-121)$\\
				\\[-1em]
				\hline
				\\[-1em]
				$0.2222$ & $0.2843$ & $0.3108$ & $0.3274$
			\end{tabular}
		\vspace{1cm}
		}
	\subcaption{Concrete values of $p_3(G_k)$ for $k=3,4,5,6$.\\ }
	\label{tab:values_p3}
\end{minipage}%
\begin{minipage}[b]{0.5\textwidth}
	\centering
	\includegraphics[width=0.75\columnwidth]{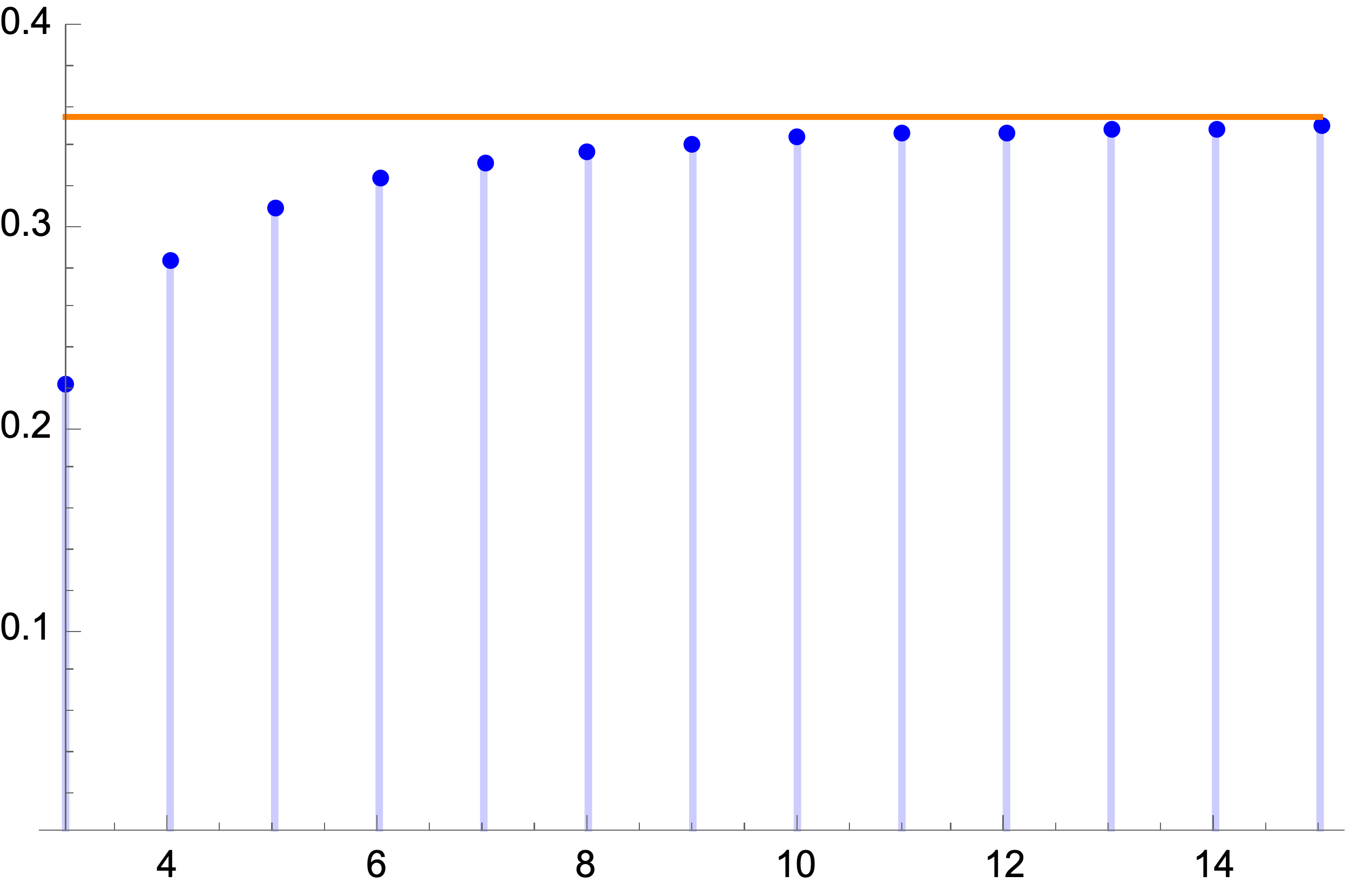}
	\subcaption{Visualization of values of $p_3(G_k)$ for $k=3,4,\dots,15$ (blue) and the limiting value $p_3(G_{\rm unif})$ (orange).}
	\label{fig:p3Gk}
\end{minipage}%
\caption{Concrete values (a) and a visualization (b) of $P_3(G_k)$ for different $k$.}
\label{fig:table_and_plot_p3}
\end{figure}

The exact and approximate values for $p_3(G_k)$ for $k\in\{3,4,5,6\}$ are summarized in the table in Figure~\ref{tab:values_p3}, some further values are visualized in Figure \ref{fig:p3Gk}. The latter also shows that, as $k\to\infty$, the value $p_3(G_k)$ tends in a monotone way to $2-{\pi^2\over 6}=p_3(G_{\rm unif})$, the triangle probability appearing in the isotropic case. This observation is confirmed in the following corollary.

\begin{corollary}\label{cor:conv}
	For $k\ge3$, let $p_3(G_k)$ be as in Theorem~\ref{thm:krichtungen}. Then $\lim\limits_{k\to\infty}p_3(G_k)= p_3(G_{\rm unif})$.
\end{corollary}

The proof of both Theorem \ref{thm:3richtungen} and Theorem \ref{thm:krichtungen} is based on the sampling procedure for the typical cell of stationary Poisson line tessellation developed in \cite{George87}. It generalizes to general directional distributions one of the stochastic constructions described in \cite{Miles73}.  To keep the paper reasonably self-contained we recall the relevant elements of this construction in the next section. Then we show in Section \ref{sec:IsotropicCase} how, using this sampling procedure, the probability $p_3(G_{\rm unif})$ can be determined. Using the same approach, the proofs of Theorem~\ref{thm:3richtungen}, Theorem~\ref{thm:krichtungen} and Corollary~\ref{cor:conv} are the content of Section~\ref{sec:proofs}. The final section of this paper provides an alternative proof of Miles' result \eqref{eq:p3Uniform} regarding the proportion of triangles in an isotropic Poisson line tessellation.

\section{Sampling random triangles}\label{sec:samplingtriangles}

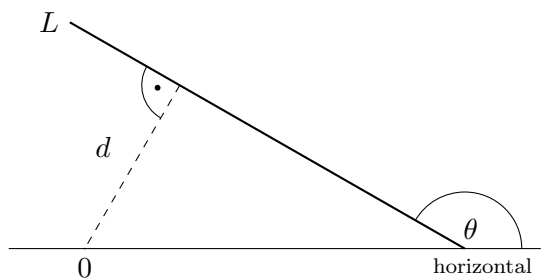
\begin{figure}[b!]
	\centering
	\begin{tikzpicture}
		\filldraw[black] (0,0) circle (0.1pt) node[below]{$0$};
		\draw[-] (-1,0)--(6,0) node[right, below, xshift=-2em]{\scriptsize horizontal};
		\draw[-, thick] (5,0)+(150:6cm)--(5,0) node[at start, left]{$L$};
		\draw[dashed] (60:2.5cm)--(0,0) node[midway, above, xshift=-1em]{$d$};
		\draw (5.75,0) arc (0:150:0.75cm) node[midway,xshift=-0.3em, yshift=-1.2em]{$\theta$};
		\draw (60:2cm) arc (240:150:0.5cm) node[midway, xshift=0.5em,yshift=0.2em]{\huge $\cdot$};
	\end{tikzpicture}
	\caption[width=0.5\columnwidth]{A line $L$ parametrized by $(d,\theta)$}
	\label{fig:parameterization_of_line}
\end{figure}

In this paper a line is parametrized by a pair $(\theta,d)$, where $d\in\RR$ is the signed distance of the line to the origin and $\theta\in[0,\pi)$ is the north-east angle this line makes with the horizontal, see Figure~\ref{fig:parameterization_of_line}. We refer to $\theta$ as the orientation angle of the line.

Following \cite{George87}, it will turn out to be convenient for us to extend the range of the possible orientation angles to the larger interval $[-\pi,\pi)$, where negative angles should be thought of modulo $\pi$. For example, we identify the orientation angles $-\pi/3$ and $2\pi/3$.

Throughout the remainder of this work, we denote random variables by a capital letter and their realizations by small ones, e.g.~$\Phi$ denotes a random angle and $\phi$ a given realization.

\subsection{General facts about Poisson line processes}

We consider a stationary Poisson line process $X=X(\gamma,G)$ with intensity $\gamma>0$ and directional distribution $G$. We assume $G$ to be non-degenerate, meaning that $G(\theta)<1$ for each $\theta\in[0,\pi)$. The following facts are taken from \cite{George87}, but see also \cite{SKM}.

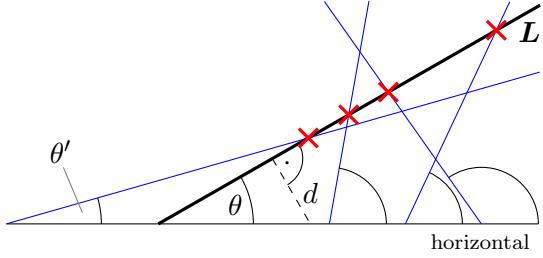
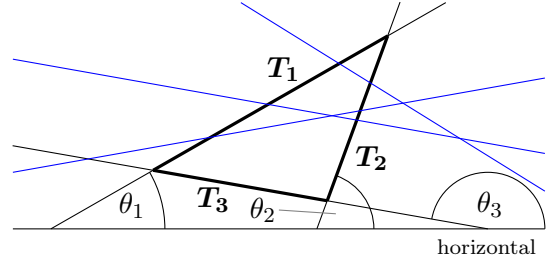
\begin{figure}[t]
	\begin{subfigure}{0.5\textwidth}
		\centering
		\begin{tikzpicture}[every mark/.append style={mark size=5pt}]
			
			\draw [-] (-2,0)--(5,0) node[right, below, xshift=-2em]{\scriptsize horizontal};
			\draw [-, line width= 1.2pt][name path=l] (0,0)--+(30:5.8cm) node[right, below, xshift=-4pt, yshift=-2pt]{$\boldsymbol L$};
			\draw (1.25,0) arc (0:30:1.25cm) node[midway, xshift=-5.5pt, yshift=-2pt]{$\theta$};
			\draw[dashed][name path=p] (2,0)+(120:1cm)--(2,0) node[right, above, xshift=0pt, yshift=4pt]{$d$};
			\path [name intersections={of=l and p, by={A}}]; 
			\draw  (A)++(30:0.4cm) arc (30:-60:0.4cm) node[midway, xshift=-0.5em,yshift=0.1em]{ $\cdot$};
			\draw [-][name path=one] [blue](-2,0)--+(16:7.3cm);
			\draw (-0.75,0) arc (0:16:1.25cm);
			\draw [-] [gray] (-1,0.15)--(-1.25,0.75) node[at end, xshift=-1pt, yshift=6pt, black]{ $\theta'$};
			\draw [-][name path=two] [blue] (3.25,0)--+(65:3.25cm);
			\draw (4,0) arc (0:65:0.75cm); 
			\draw [-][name path=three] [blue] (4.25,0)--+(125:3.6cm);
			\draw (5,0) arc (0:125:0.75cm); 
			\draw [-][name path=four] [blue] (2.25,0)--+(80:3cm);
			\draw (3,0) arc (0:80:0.75cm); 
			\path [name intersections={of=l and one, by={A}}]; 
			\draw[very thick, red] plot[mark=x] coordinates{(A)};
			\path [name intersections={of=l and two, by={A}}]; 
			\draw[very thick, red] plot[mark=x] coordinates{(A)};
			\path [name intersections={of=l and three, by={A}}]; 
			\draw[very thick, red] plot[mark=x] coordinates{(A)};
			\path [name intersections={of=l and four, by={A}}]; 
			\draw[very thick, red] plot[mark=x] coordinates{(A)};
		\end{tikzpicture}
		\subcaption{Intersection of a Poisson line process $X$ (blue) \\ with a fixed line $L=(\theta,d)$.}
		\label{fig:intersection_plp_fixed_line}
	\end{subfigure}%
	\begin{subfigure}{0.5\textwidth}
		\centering
		\begin{tikzpicture}
			
			\draw [-] (-2,0)--(5,0) node[right, below, xshift=-2em]{\scriptsize horizontal};
			\draw [-][name path=one] (-1.5,0)--+(30:6cm);
			\draw (0,0) arc (0:30:1.5cm) node[midway, xshift=-1em,yshift=-0.2em]{$\theta_1$};
			\draw [-][name path=two] (2,0)--+(70:3.2cm);
			\draw (2.75,0) arc (0:70:0.75cm);
			\draw [-] [gray] (2.25,0.2)--(1.5,0.25) node[left, xshift=0.1cm, black]{$\theta_2$};
			\draw [-][name path=four] (4.25,0)--+(170:6.35cm);
			\draw (5,0) arc (0:170:0.75cm) node[midway, xshift=-0.1em,yshift=-1.1em]{$\theta_3$};
			\path [name intersections={of=one and two, by={A}}];
			\path [name intersections={of=one and four, by={B}}]; 
			\path [name intersections={of=four and two, by={C}}]; 
			\draw[-, line width= 1.2pt] (A)--(B) node[midway, above,  xshift=5pt, yshift=4pt]{$\boldsymbol{T_1}$};
			\draw[-, line width= 1.2pt] (A)--(C) node[midway,right, xshift=-5pt, yshift=-15pt]{$\boldsymbol{T_2}$};
			\draw[-, line width= 1.2pt] (C)--(B) node[midway, below, xshift=-10pt, yshift=3.5pt]{$\boldsymbol{T_3}$};
			\draw [-][blue] (5,0.5)--(1,3);
			\draw [-][blue] (5,1)--(-2,2.25);
			\draw [-][blue] (5,2)--(-2,0.75);
		\end{tikzpicture}
		\subcaption{Lines of $X$ (blue) intersecting triangle sides $T_1$ and \\$T_2$ but not $T_3$.}
		\label{fig:intersection_plp_triangle}
	\end{subfigure}%
	\caption{Intersection of a Poisson line process with a fixed line (a) and with a triangle (b).}
\end{figure}

\paragraph{Intersection with a fixed line.} Let $L$ be a fixed line with orientation angle $\theta\in[0,\pi)$. Its intersection with $X$ is a stationary Poisson point process on $L$ with intensity $\gamma\lambda(\theta)$, where
\begin{align}\label{eq:IntersectionPP}
	\lambda(\theta) := \int_{[0,\pi)}|\sin(\theta-\theta')|\,G(\dint\theta'),
\end{align}
see Figure \ref{fig:intersection_plp_fixed_line}. Furthermore, the random orientation angles of the lines associated with these points of intersection are independent and identically distributed with common conditional density
$$
\theta'\mapsto{1\over\lambda(\theta)}|\sin(\theta-\theta')|,\qquad 0\leq\theta'<\pi,
$$
with respect to $G$.

\paragraph{Intersection of two random lines.} Let $L$ and $L'$ be two different lines from $X$, and let $(\Theta,\Theta')$ be the two orientation angles at the intersection point $L\cap L'$. Then the pair $(\Theta,\Theta')$ has joint density 
\begin{align}\label{eq:JointDensityIntersection}
	(\theta,\theta')\mapsto{1\over\lambda}|\sin(\theta-\theta')|,\qquad 0\leq\theta,\theta'<\pi,
\end{align}
with respect to the product measure $G\otimes G$ on $[0,\pi)\times[0,\pi)$, where
\begin{align}\label{eq:Lambda}
	\lambda:=\int_0^\pi\lambda(\theta)\,G(\dint\theta).
\end{align}

\paragraph{Intersection with a triangle.} Consider an arbitrary triangle $T$ in the plane with sides $T_1$, $T_2$ and $T_3$ having lengths $t_1$, $t_2$, $t_3$ and whose supporting lines have orientation angles $\theta_1$, $\theta_2$, $\theta_3$, respectively. Then the number of lines of $X$ intersecting $T$ but do not intersect $T_3$ has a Poisson distribution with mean
$$
{\gamma\over 2}(t_1\lambda(\theta_1)+t_2\lambda(\theta_2)-t_3\lambda(\theta_3)),
$$
see Figure~\ref{fig:intersection_plp_triangle}.

\subsection{Stochastic construction of a typical triangle}

A stochastic construction of the typical cell of a stationary Poisson line tessellation induced by a Poisson line process $X$ with intensity $\gamma>0$ and a general directional distribution $G$ has been described by \cite{George87} after previous works in \cite{Miles73} for the isotropic case. We rephrase it here in the special case of a triangle, i.e., we describe the distribution of the typical cell given that it is a triangle -- for brevity we refer to it as the typical triangle. Formally, the distribution $P_Z$ of the typical cell $Z$ of the Poisson line tessellation induced by $X$ is given as follows. Namely, if for a polygon $c\subset\RR^2$, $m(c)$ is the lexicographically smallest vertex, the distribution $P_Z$ of the random polygon $Z$ is given by
\begin{equation}\label{eq:TypicalCellDistributio}
P_Z(\,\cdot\,) := {1\over \EE\sum\limits_{c:m(c)\in[0,1]^2}1}\,\EE\sum\limits_{c:m(c)\in[0,1]^2}{\bf 1}\{c-m(c)\in\,\cdot\,\},
\end{equation}
where each sum runs over all cells $c$ of the Poisson line tessellation with $m(c)\in[0,1]^2$ (or any other Borel set with unit area). The distribution of the typical triangle is then the conditional distribution $P_Z(\,\cdot\,|Z\text{ is a triangle})$. 

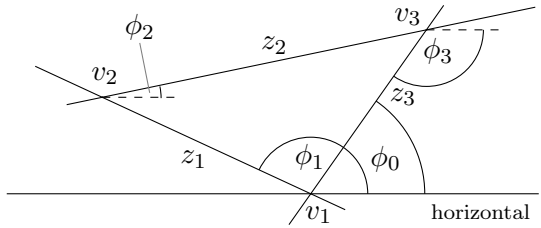
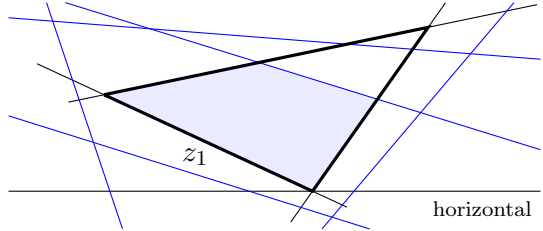
\begin{figure}[t]
	\begin{subfigure}{0.5\columnwidth}
		\centering
		\begin{tikzpicture}
			
			\draw [-] (-4,0)--(3,0) node[right, below, xshift=-2em]{\scriptsize horizontal};
			\draw[-][name path=one](0,0)+(-25:0.5cm)--(155:4cm) node[midway,left, xshift=10pt, yshift=-8pt]{ $z_1$};
			\draw[-][name path=two] (155:3cm)+(-170:0.5cm)--+(12:5.8cm) node[midway,above, xshift=-10pt, yshift=-3pt]{$z_2$};
			\draw [-][name path=three] (0,0)+(-125:0.5cm)--(55:3.05cm) node[midway,right, xshift=5pt, yshift=7pt]{$z_3$};
			\path [name intersections={of=one and two, by={A}}];
			\path [name intersections={of=two and three, by={B}}]; 
			\path [name intersections={of=one and three, by={C}}]; 
			\filldraw[black] (C) circle circle (0.1pt) node[midway, below, xshift=3pt, yshift=-1pt]{$v_1$};
			\filldraw[black] (A) circle circle (0.1pt) node[above, xshift=1pt, yshift=1pt]{$v_2$};
			\filldraw[black] (B) circle circle (0.1pt) node[left, xshift=2pt, yshift=5pt]{$v_3$};
			\draw[dashed] (A)--+(0:1);
			\draw[dashed] (B)--+(0:1);
			\draw (1.5,0) arc (0:55:1.5cm) node[midway, xshift=-0.9em,yshift=-0.6em]{$\phi_0$};
			\draw (0.75,0) arc (0:155:0.75cm) node[midway, xshift=-5pt, yshift=-7pt]{$\phi_1$};
			\draw (155:3cm)+(0.75,0) arc (0:12:0.75cm);
			\draw [-] [gray] (155:3cm)+(0.6,0.07)--(-2.25,2) node[xshift=0pt, yshift=6pt, black]{$\phi_2$};
			\draw (B)+(0.75,0) arc (0:-125:0.75cm) node[midway, xshift=-0.5em, yshift= 1em]{$\phi_3$};
		\end{tikzpicture}
		\caption{Construction of the random triangle $\Delta$.\\ \phantom{x} \\ \phantom{x}}
		\label{3 seitiges Polygon}
	\end{subfigure}%
	\begin{subfigure}{0.5\columnwidth}
		\begin{tikzpicture}
			
			\draw [-] (-4,0)--(3,0) node[right, below, xshift=-2em]{\scriptsize horizontal};
			\draw[-][name path=one](0,0)+(-25:0.5cm)--(155:4cm) node[midway,left, xshift=10pt, yshift=-8pt]{ $z_1$};
			\draw[-][name path=two] (155:3cm)+(-170:0.5cm)--+(12:5.8cm);
			\draw [-][name path=three] (0,0)+(-125:0.5cm)--(55:3.05cm);
			\draw [-][blue] [name path=b1] (-3.25,2.5)--(3,0.55);
			\draw [-][blue]  [name path=b2] (-4,1)--(0.75,-0.5);
			\draw [-][blue]  [name path=b3] (-4,2.3)--(3,1.75);
			\draw [-][blue]  [name path=b4] (0.125,-0.5)--(2.65,2.5);
			\draw [-][blue]  [name path=b5] (-3.5,2.5)--(-2.5,-0.5);
			\draw [-, line width=1.2pt] (A)--(B)--(C)--(A);
			\path [name intersections={of=one and two, by={A}}];
			\path [name intersections={of=two and three, by={B}}]; 
			\path [name intersections={of=one and three, by={C}}]; 
			\path [name intersections={of=two and b1, by={twob1}}];
			\path [name intersections={of=three and b1, by={threeb1}}];
			\path [fill=blue, fill opacity = 0.08] (A) -- (twob1) -- (threeb1) -- (C) -- (A);
		\end{tikzpicture}
		\caption{The random triangle $\Delta$ (black) and the random line process $X''$ (blue). The intersection \eqref{eq:TypicalCellDistribution} is the shaded polygon.}
		\label{fig:triangle_and_X''}
	\end{subfigure}%
	\caption{The random triangle $\Delta$ (a) and construction of the typical cell (b).}
\end{figure}

Starting with the lexicographically smallest vertex of the typical triangle, we label the vertices consecutively in clockwise direction by $v_1,v_2,v_3$. For $i\in\{1,2,3\}$ let $z_i$ be the length of the segment $\overline{v_iv_{i+1}}$, where we formally put $v_4:=v_1$. Moreover, we denote the angle between $\overline{v_iv_{i+1}}$ and the eastern horizontal at $v_i$ by $\phi_i$, see Figure \ref{3 seitiges Polygon}. Hence, $\phi_0$ denotes the initial angle. The typical triangle is completely determined by the $4$-tuple $(\Phi_0,\Phi_1,Z_1,\Phi_2)$, all other angles and edge lengths (especially $\Phi_3$, $Z_2$ and $Z_3$) can be computed from this data.

We shall now describe the (conditional) distribution of the random variables $\Phi_0,\Phi_1,Z_1$ and $\Phi_2$, which are clearly dependent. Namely,
\begin{itemize}
	\item the joint density with respect to $G\otimes G$ of $(\Phi_0,\Phi_1)$ is given by \eqref{eq:JointDensityIntersection};
	\item given $\Phi_1=\phi_1$, the intersection of $X$ with the line having orientation angle $\Phi_1$ is a stationary Poisson point process with intensity $\lambda(\phi_1)$ according to \eqref{eq:IntersectionPP}. The distance from $v_1$ to the first point of this process above the horizontal line is exponentially distributed with mean $\lambda(\phi_1)$. As a result, the conditional Lebesgue density of $Z_1$ given $\Phi_1=\phi_1$ equals 
	$$
	z_1\mapsto \lambda(\phi_1) e^{-\lambda(\phi_1)z_1},\qquad z_1>0;
	$$
	\item given $\Phi_1=\phi_1$ and $Z_1=z_1$, the random variable $\Phi_2$ has density
	$$
	\phi_2\mapsto {\sin(\phi_1-\phi_2)\over\int_{a(z_1)}^{\phi_1}\sin(\phi_1-\theta)\,G(\dint\theta)},\qquad a(z_1)\leq\phi_2<\phi_1,
	$$
	with respect to $G$. Here, $a(z_1)$ is given by
	$$
	a(z_1) := \arctan\Big({y_1\over x_1}\Big) - \pi,
	$$
	where $(x_1,y_1)$ are the coordinates of the first vertex $v_1$. 
\end{itemize}
The construction just described leads to a random triangle $\Delta$ in the plane, which is determined by the four random variables $\Phi_0,\Phi_1,\Phi_2$ and $Z_1$. It has the conditional distribution of the typical cell $Z=Z(\gamma,G)$, given that $Z$ is a triangle. To obtain from $\Delta$ the (unconditional) typical cell $Z=Z(\gamma,G)$, let $X'$ be an independent stationary Poisson line process with intensity $\gamma$ and directional distribution $G$. From $X'$ we remove all lines hitting the first edge of $\Delta$ with length $Z_1$ and call $X''$ the resulting collection of random lines, see Figure~\ref{fig:triangle_and_X''}. Then, the typical cell $Z$ has the same distribution as 
\begin{equation}\label{eq:TypicalCellDistribution}
\Delta\cap\bigcap_{L\in X''} L^+,
\end{equation}
where for each line $L$, $L^+$ denotes the closed half-space bounded by $L$ and containing the origin, see \cite{George87}.

\section{Triangle probability in the isotropic case}\label{sec:IsotropicCase}

In this section we consider the isotropic case and demonstrate how to compute $p_3:=p_3(G_{\rm unif})$ using the stochastic construction outlined in the previous section. So, let $G=G_{\rm unif}$ be the uniform distribution on $[0,\pi)$ with constant density $\theta\mapsto{\theta\over\pi}$. We also choose $\gamma=1$. It follows from \eqref{eq:IntersectionPP} and \eqref{eq:Lambda} that
$$
\lambda = \lambda(\phi) = {1\over\pi}\int_0^\pi|\sin(\phi-\theta)|\,\dint\theta = {2\over\pi},
$$
for any $\phi\in[0,\pi)$.

Due to the rotation invariance of the Poisson line tessellation in the isotropic case, the distribution of the initial angle $\Phi_0$ is irrelevant and we can just choose $\Phi_0:=0$ in the construction of the typical triangle for simplicity. Then
\begin{itemize}
	\item the random variable $\Phi_1$ has density
	$$
	\phi_1\mapsto{\pi-\phi_1\over\pi}\sin\phi_1,\qquad 0\leq\phi_1<\pi,
	$$
	which is the marginal density of the pair $(\Phi_0,\Phi_1)$ with respect to the second coordinate;
	\item the random variable $Z_1$ is independent of $\Phi_1$ and has density
	$$
	z_1\mapsto {2\over\pi}e^{-{2\over\pi}z_1},\qquad z_1>0;
	$$
	\item the random variable $\Phi_2$ only depends on $\Phi_1$ and, given $\Phi_1=\phi_1$, has conditional density
	$$
	\phi_2\mapsto{1\over 2}\sin(\phi_1-\phi_2),\qquad \phi_1-\pi\leq\phi_2 < 0,
	$$ 
	since $x_1=z_1\cos\phi_1$, $y_1=z_1\sin\phi_1$, which in turn implies that $a(z_1)=\phi_1-\pi$, independently of $z_1$.
\end{itemize}
Given these distributions, the probability $p_3$ that the typical cell is a triangle can now be written as follows:
\begin{align*}
	p_3 &= \int_0^\pi\int_0^\infty\int_0^{\phi_1-\pi} e^{-{1\over 2}(\lambda(\phi_2)z_2+\lambda(\phi_3)z_3-\lambda(\phi_1)z_1)}\\
	&\qquad\qquad\qquad\times {\pi-\phi_1\over\pi}\sin\phi_1\times {2\over\pi}e^{-{2\over\pi}z_1}\times {1\over 2}\sin(\phi_1-\phi_2)\,\dint\phi_2\dint z_1\dint\phi_1.
\end{align*}
In fact, in order to ensure that the typical cell is a triangle, we need to ensure that after the stochastic construction of the typical triangle, giving $\Phi_1=\phi_1$, $Z_1=z_1$ and $\Phi_2=\phi_2$, the two edges with length $z_2$ and $z_3$ are not intersected by lines of the random line process $X''$, recall \eqref{eq:TypicalCellDistribution}. Thus, by the intersection-with-a-triangle-property the above event has probability $\exp(-{1\over 2}(\lambda(\phi_2)z_2+\lambda(\phi_3)z_3-\lambda(\phi_1)z_1))$, which is the probability that a Poisson random variable with mean ${1\over 2}(\lambda(\phi_2)z_2+\lambda(\phi_3)z_3-\lambda(\phi_1)z_1)$ takes the value zero. The other terms in the above integral representation are just the densities of the random variables $\Phi_1$, $Z_1$ and $\Phi_2$.

It is not difficult to verify that
\begin{align}
	z_2 &= -z_1{\sin\phi_1\over\sin\phi_2},\label{eq:z2}\\
	z_3 &= z_2\cos\phi_1-z_1{\sin\phi_1\over\sin\phi_2}\cos\phi_2,\label{eq:z3}\\
	\phi_3 &= \phi_0-\pi.\label{eq:phi3}
\end{align}
This yields
$$
z_1+z_2+z_3 = z_1\,{\sin\phi_2-\sin\phi_1-\sin(\phi_1-\phi_2)\over\sin\phi_2}.
$$
Inserting this together with the values of $\lambda(\phi_1)=\lambda(\phi_2)=\lambda(\phi_3)={2\over\pi}$ we arrive at
\begin{align*}
	p_3 &=\int_0^{\pi}{\pi-\phi_1\over\pi^2}\int_{\phi_1-\pi}^0\sin\phi_1\sin(\phi_1-\phi_2)\int_0^\infty e^{-{z_1\over\pi}{\sin\phi_2-\sin\phi_1-\sin(\phi_1-\phi_2)\over\sin\phi_2}}\,\dint z_1\dint\phi_2\dint\phi_1.
\end{align*}
Solving the innermost integral leads, after simplification of the resulting expression, to
\begin{align}\label{eq:IsoInt2}
		p_3 &=\int_0^{\pi}{\pi-\phi_1\over\pi}\int_{\phi_1-\pi}^0{\sin\phi_1\sin\phi_2\sin(\phi_1-\phi_2)\over\sin\phi_2-\sin\phi_1-\sin(\phi_1-\phi_2)}\,\dint\phi_2\dint\phi_1.
\end{align}
Rewriting the integrand by means of trigonometric identities gives
\begin{equation}\label{eq:IsotropicIntegral}
\begin{split}
	p_3 &= \int_0^\pi{\pi-\phi_1\over\pi}\int_0^{\pi-\phi_1}\sin{\phi_1\over 2}\Big(\sin\Big(\phi_2+{\phi_1\over 2}\Big)-\sin{\phi_1\over 2}\Big)\,\dint\phi_2\dint\phi_1\\
	&={1\over 2\pi}\int_0^\pi(\pi-\phi_1)(2\sin\phi_1-(\pi-\phi_1)(1-\cos\phi_1))\,\dint\phi_1= 2 - {\pi^2\over 6},
\end{split}
\end{equation}
all details of the computation were carried out in \cite{Krecklenberg}.

\section{Proofs}\label{sec:proofs}

\subsection{Triangle probability in the $\boldsymbol{G_3(p,q)}$-case: Proof of Theorem \ref{thm:3richtungen}}\label{sec:proof3rich}

In this section we compute the triangle probability $p_3:=p_3(G_3(p,q))$ if the underlying directional distribution is given by \eqref{eq:G3pq}, again using the stochastic construction of the typical cell. We recall that we choose $\gamma=1$ as our intensity. 

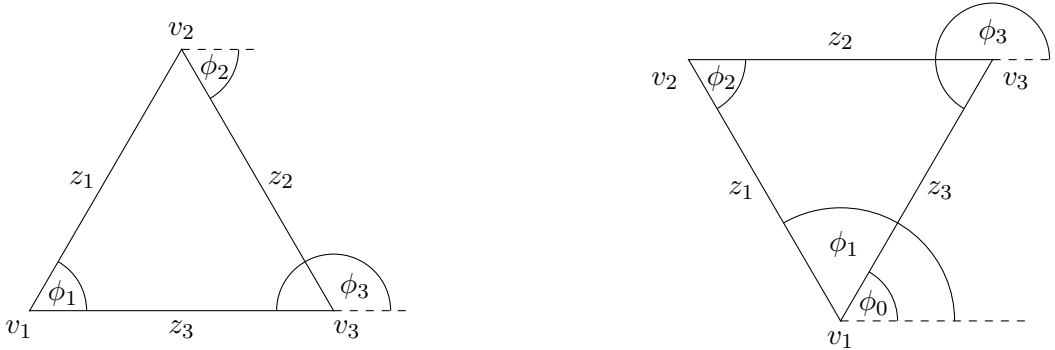
\begin{figure}
	\begin{subfigure}{0.5\textwidth}
		\centering
		\begin{tikzpicture}
			
			\draw[-](0,0)--(60:4cm) node[midway ,left]{$z_1$};
			\draw[-](4,0)--(60:4cm) node[midway, right]{$z_2$};
			\draw[-](0,0)--(4,0) node[midway,below]{$z_3$}; 
			\draw[dashed](60:4cm)--+(1,0);
			\draw [dashed] (4,0)--(5,0);
			\draw (0.75,0) arc (0:60:0.75cm) node[midway, xshift=-6pt, yshift=-4pt]{$\phi_1$};
			\draw (60:4cm)+(0.75,0) arc (360:300:0.75cm)node[midway, xshift=-6pt, yshift=4pt]{$\phi_2$};
			\draw (4.75,0) arc (0:180:0.75cm)node[midway, xshift=8pt, yshift=-12pt]{$\phi_3$};
			\filldraw[black] (0,0) circle (0.1pt) node[midway, below, xshift=-4pt]{$v_1$};
			\filldraw[black] (60:4cm) circle (0.1pt) node[above]{$v_2$};
			\filldraw[black] (4,0) circle (0.1pt) node[below, xshift=5pt]{$v_3$};
		\end{tikzpicture}
		\caption{With angles $\phi_0=0$, $\phi_1=\frac{1}{3}\pi$ and $\phi_2=-\frac{1}{3}\pi$.}
	\end{subfigure}
	\begin{subfigure}{0.5\textwidth}
		\centering
		\begin{tikzpicture}
			
			\draw[-](120:4cm)--(0,0) node[midway,left]{$z_1$};
			\draw[-](120:4cm)--+(4,0) node[midway, above]{$z_2$};
			\draw[-](60:4cm)--(0,0)node[midway, right]{$z_3$};
			\draw[dashed](0,0)--(2.5,0);
			\draw [dashed] (60:4cm)--+(1,0);
			\draw (0.75,0) arc (0:60:0.75cm) node[midway, xshift=-6pt, yshift=-4pt]{$\phi_0$};
			\draw (1.5,0) arc (0:120:1.5cm) node[midway, xshift=-20pt, yshift=-8pt]{$\phi_1$};
			\draw (120:4cm)+(0.75,0) arc (360:300:0.75cm)node[midway, xshift=-6pt, yshift=4pt]{$\phi_2$};
			\draw (60:4cm)+(0.75,0) arc (0:240:0.75cm) node[midway, xshift=11pt, yshift=-8pt]{$\phi_3$};
			\filldraw[black] (0,0) circle (0.1pt) node[midway, below]{$v_1$};
			\filldraw[black] (120:4cm) circle (0.1pt) node[left, yshift=-8pt]{$v_2$};
			\filldraw[black] (60:4cm) circle (0.1pt) node[right, yshift=-8pt]{$v_3$};
		\end{tikzpicture}
		\caption{With angles $\phi_0=\frac{1}{3}\pi$, $\phi_1=\frac{2}{3}\pi$ and $\phi_2=0$.}
	\end{subfigure}
	\caption{The two possible triangles in a Poisson line tessellation with directional distribution $G_3(p,q)$.}
	\label{fig:2DreieckeG3}
\end{figure}

Before we actually compute $p_3$, we deal with the possible constructions for triangles with only three edge directions corresponding to the orientation angles $0$, $\pi\over 3$ and $2\pi\over 3$. In fact, we only have two ways to construct a triangle with these orientation angles as demonstrated in Figure \ref{fig:2DreieckeG3}. Further, writing $G$ for $G_3(p,q)$ for brevity, we can now compute
\begin{align*}
	\lambda(0) &= \int_0^\pi|\sin\theta|\,G(\dint\theta) = q\sin{\pi\over 3} + (1-p-q)\sin{2\pi\over 3} = {\sqrt{3}\over 2}(1-p),\\
	\lambda\Big({\pi\over 3}\Big) &=  \int_0^\pi\Big|\sin\Big(\theta-{\pi\over 3}\Big)\Big|\,G(\dint\theta) 
	=p\Big|\sin\Big(-{\pi\over 3}\Big)\Big|+(1-p-q)\Big|\sin\Big({2\pi\over 3}-{\pi\over 3}\Big)\Big|
	={\sqrt{3}\over 2}(1-q),\\
	\lambda\Big({2\pi\over 3}\Big) &=  \int_0^\pi\Big|\sin\Big(\theta-{2\pi\over 3}\Big)\Big|\,G(\dint\theta) 
	=p\Big|\sin\Big(-{2\pi\over 3}\Big)\Big|+q\Big|\sin\Big({\pi\over 3}-{2\pi\over 3}\Big)\Big|
	={\sqrt{3}\over 2}(p+q),
\end{align*}
according to \eqref{eq:IntersectionPP}, which implies that
\begin{align*}
	\lambda &= \int_0^\pi\lambda(\theta)\,G(\dint\theta)=p\lambda(0)+q\lambda\Big({\pi\over 3}\Big)+(1-p-q)\lambda\Big({2\pi\over 3}\Big)=\sqrt{3}(p+q-p^2-q^2-pq).
\end{align*}
Moreover, from \eqref{eq:JointDensityIntersection} it follows that the pair $(\Phi_0,\Phi_1)$ has joint density
\begin{align*}
	(\phi_0,\phi_1)\mapsto &{2\over \sqrt{3}(p+q-p^2-q^2-pq)}\sin(\phi_0-\phi_1),\qquad 0\leq\phi_0,\phi_1<\pi,
\end{align*}
with respect to $G\otimes G$. Given $\Phi_1=\phi_1$, the random variable $Z_1$ is exponentially distributed with mean $\lambda(\phi_1)$. Finally, as in the isotropic case, we have $a(z_1)=\phi_1-\pi$ and so the random variable $\Phi_2$ has conditional density
$$
\phi_2\mapsto{1\over\lambda(\phi_1)}\sin(\phi_1-\phi_2),\qquad \phi_1-\pi\leq\phi_2<0,
$$
with respect to $G$, given $\Phi_1=\phi_1$.

With the same argument as in the isotropic case, we can now represent the triangle probability as follows:
\begin{align*}
		p_3 &= \int_0^\pi\int_0^\pi\int_{\phi_1-\pi}^0\int_0^\infty e^{-{1\over 2}(\lambda(\phi_2)z_2+\lambda(\phi_3)z_3-\lambda(\phi_1)z_1)}\times{2\over \sqrt{3}(p+q-p^2-q^2-pq)}\sin(\phi_0-\phi_1)\\
		&\qquad\times{1\over\lambda(\phi_1)}\sin(\phi_1-\phi_2)\times\lambda(\phi_1)e^{-\lambda(\phi_1)z_1}\,\dint z_1G(\dint\phi_2)(G\otimes G)(\dint(\phi_0,\phi_1))\\
		&= \int_0^\pi\int_0^\pi\int_{\phi_1-\pi}^0\int_0^\infty e^{-{1\over 2}(\lambda(\phi_1)z_1+\lambda(\phi_2)z_2+\lambda(\phi_3)z_3)}\,\times{2\over \sqrt{3}(p+q-p^2-q^2-pq)}\\
		&\qquad\times\sin(\phi_0-\phi_1)\sin(\phi_1-\phi_2)\,\dint z_1G(\dint\phi_2)(G\otimes G)(\dint(\phi_0,\phi_1));
\end{align*}
the term $e^{-{1\over 2}(\lambda(\phi_2)z_2+\lambda(\phi_3)z_3-\lambda(\phi_1)z_1)}$ represents the probability that after the stochastic construction of the typical triangle the random line process $X''$ does not intersect the two edges with lengths $z_2$ and $z_2$, whereas the other terms are the (conditional) densities of $(\Phi_0,\Phi_1)$, $\Phi_2$ and $Z_1$. From the discussion at the beginning of this section we know the three outer integrals are just a sum of two terms corresponding to the following angles:
\begin{alignat*}{5}
	&\text{Case 1:}\quad&&\phi_0 =0,\quad&&\phi_1={\pi\over 3},\quad&&\phi_2 ={2\pi\over 3},\quad&&\phi_3 =-\pi,\\
	&\text{Case 2:}\quad &&\phi_0={\pi\over 3},&&\phi_1={2\pi\over 3},&&\phi_2 ={0},&&\phi_3 =-{2\pi\over 3}.
\end{alignat*}
In both cases, using \eqref{eq:z2} and \eqref{eq:z3}, we conclude that, given $\Phi_1=\phi_1$, $Z_1=z_1$ and $\Phi_2=\phi_2$, we have $z_1=z_2=z_3$, formally confirming that we are dealing with regular triangles. Moreover, in both cases we have $\lambda(\phi_1)+\lambda(\phi_2)+\lambda(\phi_3) = \sqrt{3}$, implying that
\begin{align*}
		&\int_0^\infty e^{-{1\over 2}(\lambda(\phi_1)z_1+\lambda(\phi_2)z_2+\lambda(\phi_3)z_3)}\,\dint z_1=\int_0^\infty e^{-{\sqrt{3}\over 2}z_1}\,\dint z_1 = {2\over\sqrt{3}}.
\end{align*}
Hence,
\begin{align*}
		p_3 = {2\over\sqrt{3}}&\times{2\over \sqrt{3}(p+q-p^2-q^2-pq)}\\
		&\times\int_0^\pi\int_0^\pi\int_{\phi_1-\pi}^0\sin(\phi_0-\phi_1)\sin(\phi_1-\phi_2)\,G(\dint\phi_2)(G\otimes G)(\dint(\phi_0,\phi_1)).
\end{align*}
Finally, in case 1, which has weight $pq(1-p-q)$, the integrand equals $3/4$, and in case 2, which has weight $p(1-p-q)p$, the integrand equals $3/4$ as well. This eventually leads to
\begin{align*}
		p_3 &= {2\over\sqrt{3}}\times{2\over \sqrt{3}(p+q-p^2-q^2-pq)}\times 2\times pq(1-p-q)\times{3\over 4}={2pq(1-p-q)\over p+q-p^2-q^2-pq}
\end{align*}
and concludes the proof of the first part of Theorem \ref{thm:3richtungen}.

For the second part, define the function
$$
F(p,q) := {2pq(1-p-q)\over p+q-p^2-q^2-pq}
$$
on the domain $D:=\{(p,q)\in(0,1)^2:0<p+q<1\}$ whose gradient is
\begin{align*}
	&{\rm grad}F(p,q) = {2\over (p+q-p^2-q^2-pq)^2}\big(q(1-p-q)-pq,p(1-p-q)-pq\big).
\end{align*}
Solving ${\rm grad}F(p,q)=(0,0)$ leads to the only solution $(p,q)=(1/3,1/3)$ on $D$. One can easily check that this is indeed the global maximum of $F(p,q)$ on $D$. Since $p_3(G_3(1/3,1/3))=2/9$, the proof of Theorem \ref{thm:3richtungen} is complete.\qed

\subsection{Triangle probability in the case of $\boldsymbol{k}$ directions: Proof of Theorem \ref{thm:krichtungen}}\label{sec:proofkrich}

Recall the construction of a typical triangle based on the random angles $\Phi_0,\Phi_1,\Phi_2$ and the random edge length $Z_1$. Since the Poisson line tessellation with directional distribution $G_k$ is $G_k$-pseudo isotropic, the initial angle $\Phi_0$ is irrelevant and we can just take $\Phi_0=0$. Moreover, recall that $\Phi_3=\Phi_0-\pi=-\pi$. It is now a crucial observation that the stochastic construction described above leads to a triangle if and only if 
$$
(\phi_1,\phi_2)\in\Big\{\Big({i\pi\over k},-{j\pi\over k}\Big):{1\leq i\leq k-2\atop 1\leq j\leq k-i-1}\Big\},
$$
since the angle sum of a triangle is equal to $\pi$ and since we require the vertex $v_1$ to be the lexicographically smallest vertex of the triangle. Moreover, for fixed $1\leq i\leq k-2$ each such triangle can be rotated by the angles $0,{\pi\over k},\ldots,{(k-i-1)\pi\over k}$ to yield another admissible triangle.

We determine now the distribution of the relevant random variables and start with $\Phi_1$. According to \eqref{eq:IntersectionPP} and using the identity for sums of sines in arithmetic progressions from \cite{Knapp} (with $a=0$ and $d=\pi/k$ there) we have
\begin{align*}
	\lambda(0) &= \int_0^\pi|\sin(\theta)|\,G_k(\dint\theta) = {1\over k}\sum_{\ell=0}^{k-1}\sin{\ell\pi\over k} = {1\over k}{\sin{(k-1)\pi\over 2k}\over\sin{\pi\over 2k}} = {1\over k}\cot{\pi\over 2k},
\end{align*}
and because of $G_k$-pseudo isotropy we also have $\lambda({\pi\over k})=\ldots=\lambda({(k-1)\pi\over k})=\lambda(0)$. Thus, from \eqref{eq:Lambda} it follows that
\begin{align*}
	\lambda = \int_0^\pi\lambda(\theta)\,G_k(\dint\theta) = {1\over k}\cot{\pi\over 2k}.
\end{align*}
We can now conclude from \eqref{eq:JointDensityIntersection} that the pair $(\Phi_0,\Phi_1)$ has joint density
\begin{align*}
	(\phi_0,\phi_1)\mapsto {2k\over\cot{\pi\over 2k}}\sin(\phi_1-\phi_0),\qquad 0\leq\phi_0,\phi_1<\pi,
\end{align*}
with respect to $G_k\otimes G_k$. Integration with respect to $\phi_0$ yields now the marginal density
\begin{align*}
	\phi_1\mapsto &{2k\over\cot{\pi\over 2k}}\int_0^{\pi-\phi_1} \sin(\phi_1-\phi_0)\,G_k(\dint\phi_0)={2\Sigma_k(\phi_1)\over \cot{\pi\over 2k}}\sin(\phi_1),\qquad 0\leq\phi_1<\pi,
\end{align*}
of $\Phi_1$ with respect to $G_k$, where
\begin{align}\label{eq:sigma_k}
\Sigma_k(\phi_1) := \sum_{\phi_0\in\{0,{\pi\over k},\ldots,{(k-1)\pi\over k}\}\atop \phi_0<\pi-\phi_1}1.
\end{align}
The distribution of $Z_1$ is an exponential distribution with mean ${1\over k}\cot{\pi\over 2k}$ and so $Z_1$ has density
$$
z_1 \mapsto {1\over k}\cot{\pi\over 2k} e^{-{1\over k}\cot{\pi\over 2k}\,z_1},\qquad z_1>0,
$$
with respect to the Lebesgue measure. Finally, we deal with the conditional distribution of $\Phi_2$ given $\Phi_1$. As above, we have that the conditional density with respect to $G_k$ of $\Phi_2$ given $\Phi_1=\phi_1$ equals
$$
\phi_2\mapsto {\sin(\phi_1-\phi_2)\over\int_{\phi_1-\pi}^{\phi_1}\sin(\phi_1-\phi)\,G_k(\dint\phi)},\qquad \phi_1-\pi\leq\phi_2<\phi_1.
$$
Since the integral in the denominator is just $\lambda(\phi_1)$, we arrive at
$$
\phi_2 \mapsto {k\over\cot{\pi\over 2k}}\,\sin(\phi_1-\phi_2),\qquad \phi_1-\pi\leq\phi_2<\phi_1,
$$
for the conditional density of $\Phi_2$.

As in the two previous sections, we can now express $p_3:=p_3(G_k)$ as follows:
\begin{align*}
	p_3 &= \int_0^\pi\int_{\phi_1-\pi}^0\int_0^\infty {2\Sigma_k(\phi_1)\over\cot{\pi\over 2k}}\sin(\phi_1)\times {1\over k}\cot{\pi\over 2k} e^{-{1\over k}\cot{\pi\over 2k}\,z_1}\times {k\over\cot{\pi\over 2k}}\,\sin(\phi_1-\phi_2)\\
	&\hspace{5cm}\times e^{-{1\over 2}(\lambda(\phi_2)z_2+\lambda(\phi_3)z_3-\lambda(\phi_1)z_1)}\,\dint z_1G_k(\dint\phi_2)G_k(\dint\phi_1)\\
	&= {2\over\cot{\pi\over 2k}}\int_0^\pi\int_{\phi_1-\pi}^0 \Sigma_k(\phi_1)\sin(\phi_1)\sin(\phi_1-\phi_2)\\
	&\hspace{5cm}\times\int_0^\infty e^{-{1\over 2k}\cot{\pi\over 2k}(z_1+z_2+z_3)}\,\dint z_1G_k(\dint\phi_2)G_k(\dint\phi_1),
\end{align*}
where in the last step we used that $\lambda(\phi)=\lambda(0)$ for all angles $\phi$ in the support of $G_k$. 

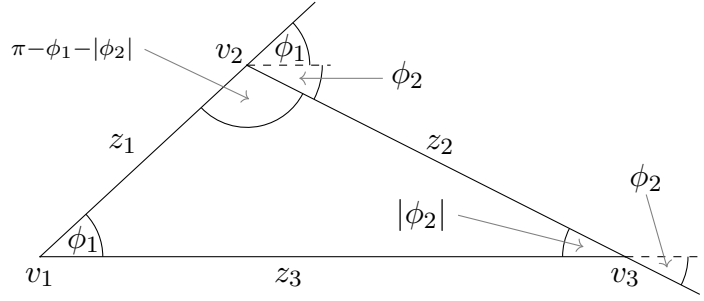
\begin{figure}[t]
	\centering
		\begin{tikzpicture}[scale=1.1, transform shape]
			
			\draw[-] (0,0)--+(43:4.5cm) node[midway, left, xshift=-10, yshift=-5]{$z_1$};
			\draw[-] (7,0)++(-27:1cm)--+(153:6.085cm) node[midway, right, xshift=-20pt, yshift=12pt]{$z_2$};
			\draw[-] (0,0)--(7,0) node[midway, below, xshift=-15pt]{$z_3$};
			\draw[dashed](7,0)--+(1,0);
			\draw[dashed](43:3.385cm)--+(1,0);
			\draw (0.75,0) arc (0:43:0.75cm) node [midway, xshift=-5pt, yshift=-2pt]{$\phi_1$}; 
			\draw (43:3.385cm)++(0.9,0) arc (0:-27:0.9cm); 
			\draw [<-, gray] (3.1,2.15)--(4.1,2.15) node[right, xshift=-1pt, yshift=1pt, black]{$\phi_2$};
			\filldraw (0,0) circle (0.1pt) node[below]{$v_1$};
			\filldraw (43:3.385cm) circle (0.1pt) node[left, xshift=3pt, yshift=4pt]{$v_2$};
			\filldraw (7,0) circle (0.1pt) node[below]{$v_3$};
			\draw (43:4.135cm) arc (43:0:0.75cm) node[midway, xshift=-5pt, yshift=-2pt]{$\phi_1$};
			\draw (7,0)++(-27:0.75cm) arc (-27:0:0.75cm);
			\draw [<-, gray] (7.5,-0.13)--(7.25,0.75) node[at end, above, yshift=-3pt, black]{$\phi_2$};
			\draw(43:2.635cm) arc (-137:-27:0.75cm);
			\draw [<-, gray] (2.5,1.9)--(1.25,2.5) node[left, black]{$\substack{\pi-\phi_1-\vert\phi_2\vert}$};
			\draw (7,0)++(153:0.75cm) arc (153:180:0.75cm);
			\draw [<-, gray] (6.5,0.12)--(5,0.5) node[at end, left, black]{$\vert\phi_2\vert$};
		\end{tikzpicture}
		\caption{Determination of $z_2$ and $z_3$.}
		\label{fig:z2z3}
		
\end{figure}

To determine $z_2$ and $z_3$ we can use the law of sines as illustrated in Figure \ref{fig:z2z3}.
If $\phi_1={i\pi\over k}$, $1\leq i\leq n-2$, and $\phi_2=-{j\pi\over k}$, $1\leq j\leq n-i-1$, this yields
\begin{align*}
	z_2 = z_1{\sin{i\pi\over k}\over\sin{j\pi\over k}}\qquad\text{and}\qquad
	z_3 = z_1{\sin{{(i+j)}\pi\over k}\over\sin{j\pi\over k}}.
\end{align*}
Thus,
\begin{align*}
	&{1\over 2k}\cot{\pi\over 2k}(z_1+z_2+z_3) = {1\over 2k}\cot{\pi\over 2k}{\sin{i\pi\over k}+\sin{j\pi\over k}+\sin{(i+j)\pi\over k}\over\sin{j\pi\over k}}z_1,
\end{align*}
and the integral with respect to $z_1$ evaluates to
$$
\int_0^\infty e^{-{1\over 2k}\cot{\pi\over 2k}(z_1+z_2+z_3)}\,\dint z_1={2k\sin{j\pi\over k}\over \cot{\pi\over 2k}\big(\sin{i\pi\over k}+\sin{j\pi\over k}+\sin{(i+j)\pi\over k}\big)}.
$$
Plugging this back into the expression for $p_3$, we see that
\begin{align}\label{eq:p3GkRepresentation}
	p_3 &= {4k\over\cot^2{\pi\over 2k}}{1\over k^2}\sum_{i=1}^{k-2}\Sigma_k\Big({i\pi\over k}\Big)\sum_{j=1}^{k-i-1}{\sin{i\pi\over k}\,\sin{j\pi\over k}\,\sin{(i+j)\pi\over k}\over \sin{i\pi\over k}+\sin{j\pi\over k}+\sin{(i+j)\pi\over k}}.
\end{align}
Using that $\Sigma_k\big({i\pi\over k}\big)=k-i$, we can complete the proof of Theorem \ref{thm:krichtungen}.\qed

\subsection{The convergence to the isotropic case: Proof of Corollary~\ref{cor:conv}}\label{sec:proof_conv_to_unif}

We start with the observation that 
\begin{align}\label{eq:help_proof_conv_1}
	\frac{4}{k}\tan^2\frac{\pi}{2k}=\frac{\pi^2}{k^3}+O(k^{-5}),
\end{align}
as $k\to\infty$. Combining this with the representation for $p_3(G_k)$ in Theorem \ref{thm:krichtungen} implies
\begin{align*}
	\lim_{k\to\infty} p_3(G_k) &= \lim_{k\to\infty} {4\over k}\tan^2{\pi\over 2k}\sum_{i=1}^{k-2}\Bigg[(k-i)\sum_{j=1}^{k-i-1}{\sin{i\pi\over k}\,\sin{j\pi\over k}\,\sin{(i+j)\pi\over k}\over\sin{i\pi\over k}+\sin{j\pi\over k}+\sin{(i+j)\pi\over k}}\Bigg]\\
	&= \lim_{k\to\infty} \frac{\pi^2}{k^3}\sum_{i=1}^{k-2}\Bigg[(k-i)\sum_{j=1}^{k-i-1}{\sin{i\pi\over k}\,\sin{j\pi\over k}\,\sin{(i+j)\pi\over k}\over\sin{i\pi\over k}+\sin{j\pi\over k}+\sin{(i+j)\pi\over k}}\Bigg]\\
	&= \lim_{k\to\infty} \pi^2\frac{1}{k}\sum_{i=1}^{k-2}\Bigg[\Big(1-\frac{i}{k}\Big)\frac{1}{k}\sum_{j=1}^{k-i-1}{\sin{i\pi\over k}\,\sin{j\pi\over k}\,\sin{(i+j)\pi\over k}\over\sin{i\pi\over k}+\sin{j\pi\over k}+\sin{(i+j)\pi\over k}}\Bigg].
\end{align*}
Interpreting the two sums as Riemann sums with $\frac{i}{k}\to \dint t$ and $\frac{j}{k}\to \dint s$, as $k\to\infty$, and noting that the condition $j\leq k-i-1$ asymptotically translates to $s< 1-t$, we conclude that
\begin{align}\label{eq:Limp3Gk}
	\lim_{k\to\infty} p_3(G_k)=\pi^2\int_0^1 (1-t)\int_0^{1-t}{\sin(\pi t)\,\sin(\pi s)\,\sin((t+s)\pi)\over\sin(\pi t)+\sin(\pi s)+\sin((t+s)\pi)}\,\dint s\dint t.
\end{align}
This, up to the substitutions $u=\pi t$ and $v=-\pi s$, is exactly the integral expression for $p_3(G_{\rm unif})$ we encountered already in \eqref{eq:IsoInt2}. This completes the argument. \qed

\makeatletter
\renewcommand*\maketag@@@[1]{\hbox{\m@th #1}}
\makeatother

\section{Alternative proof of Miles' result {\eqref{eq:p3Uniform}}}\label{sec:alt_miles}

As mentioned in the introduction, it is known from \cite{Miles64} that $p_3(G_{\rm unif})=2-\frac{\pi^2}{6}$. In this section, we use our Theorem~\ref{thm:krichtungen} to give a `continuous-mapping-type' argument leading to the same result. Our strategy is to prove that the weak convergence of $G_k$ to $G_{\rm unif}$ implies the convergence of $p_3(G_k)$ to $p_3(G_{\rm unif})$, as $k\to\infty$. To conclude, we can then use Corollary \ref{cor:conv}, which shows that $p_3(G_{\rm unif})=\lim\limits_{k\to\infty}p_3(G_k)$. The value of this limit is given by the integral \eqref{eq:Limp3Gk}, which we evaluated to $2-{\pi^2\over 6}$ in \eqref{eq:IsotropicIntegral}. The approach can be summarized in the following chain of equalities, in which $\lim\limits^w$ stands for the weak limit of probability measures:
$$
p_3(\lim\limits_{k\to\infty}^wG_k) \overset{{\rm shown\ below}}{=}  \lim_{k\to\infty}p_3(G_k) \overset{{\rm Corollary\ \ref{cor:conv}}}{=} p_3(G_{\rm unif})   \overset{\eqref{eq:IsotropicIntegral}}{=} 2-{\pi^2\over 6}.
$$
To prove the first equality, we recall that the weak convergence of $G_k$ to $G_{\rm unif}$ implies the weak convergence of the product measures $G_k\otimes G_k\otimes G_k$ to $G_{\rm unif}\otimes G_{\rm unif}\otimes G_{\rm unif}$, see \cite[Proposition 2.7.7]{Boga}. For each $k\geq 3$ the triangle probability $p_3(G_k)$ can be represented as the integral
$$
p_3(G_k) = \int_{[0,\pi)\times[0,\pi)\times[0,\pi)} f_k(\phi_0,\phi_1,\phi_2)\,(G_k\otimes G_k\otimes G_k)(\dint(\phi_0,\phi_1,\phi_2))
$$
with the function $f_k:[0,\pi)\times[0,\pi)\times[0,\pi)\to\RR$ given by
\begin{align}\label{eq:Functionfk}
f_k(\phi_0,\phi_1,\phi_2):=4k^2\tan^2{\pi\over2k}T(\phi_0,\phi_1,\phi_2){\bf 1}\{\phi_0<\pi-\phi_1,\phi_1<\phi_2\},
\end{align}
where
$$
T(\phi_0,\phi_1,\phi_2):={\sin(\phi_0-\phi_1)\sin\phi_1\,\sin\phi_2\,\sin(\phi_1-\phi_2) \over \sin\phi_1+\sin\phi_2+\sin(\phi_1-\phi_2)}.
$$
Note that if the integration with respect to $\phi_0$ is carried out, we precisely arrive at \eqref{eq:p3GkRepresentation}. Since $4k^2\tan^2{\pi\over2k}\to \pi^2$ as $k\to\infty$ by \eqref{eq:help_proof_conv_1}, we have that
\begin{align}\label{eq:19922A}
	\lim_{k\to\infty}f_k(\phi_0,\phi_1,\phi_2) = \pi^2 T(\phi_0,\phi_1,\phi_2){\bf 1}\{\phi_0<\pi-\phi_1,\phi_1<\phi_2\}=:f(\phi_0,\phi_1,\phi_2)
\end{align}
pointwise on $[0,\pi)\times[0,\pi)\times[0,\pi)$. To complete the proof, it remains to verify that
\begin{align*}
&\lim_{k\to\infty} \int_{[0,\pi)\times[0,\pi)\times[0,\pi)} f_k(\phi_0,\phi_1,\phi_2)\,(G_k\otimes G_k\otimes G_k)(\dint(\phi_0,\phi_1,\phi_2))\\
&\qquad =  \int_{[0,\pi)\times[0,\pi)\times[0,\pi)} \lim_{k\to\infty} f_k(\phi_0,\phi_1,\phi_2)\;\lim\limits_{k\to\infty}^w(G_k\otimes G_k\otimes G_k)(\dint(\phi_0,\phi_1,\phi_2)) \\
& \qquad= \int_{[0,\pi)\times[0,\pi)\times[0,\pi)} f(\phi_0,\phi_1,\phi_2)\,(G_{\rm unif}\otimes G_{\rm unif}\otimes G_{\rm unif})(\dint(\phi_0,\phi_1,\phi_2)).
\end{align*}
For this, since $4k^2\tan^2{\pi\over2k}\to \pi^2$, we only need to prove that
\begin{align*}
	&\lim_{k\to\infty} \int_{[0,\pi)\times[0,\pi)\times[0,\pi)} T(\phi_0,\phi_1,\phi_2){\bf 1}\{\phi_0<\pi-\phi_1,\phi_1<\phi_2\}\,(G_k\otimes G_k\otimes G_k)(\dint(\phi_0,\phi_1,\phi_2))\\
	& \qquad= \int_{[0,\pi)\times[0,\pi)\times[0,\pi)} T(\phi_0,\phi_1,\phi_2){\bf 1}\{\phi_0<\pi-\phi_1,\phi_1<\phi_2\}\,(G_{\rm unif}\otimes G_{\rm unif}\otimes G_{\rm unif})(\dint(\phi_0,\phi_1,\phi_2)).
\end{align*}
However, since the function $(\phi_0,\phi_1,\phi_2)\mapsto T(\phi_0,\phi_1,\phi_2){\bf 1}\{\phi_0<\pi-\phi_1,\phi_1<\phi_2\}$ is bounded and $G_{\rm unif}\otimes G_{\rm unif}\otimes G_{\rm unif}$-almost everywhere continuous on $[0,\pi)\times[0,\pi)\times[0,\pi)$, this follows from \cite[Corollary 2.2.10]{Boga} and the argument is complete.

\subsection*{Acknowledgment}
We are grateful to Tom Kaufmann and Daniel Rosen for inspiring ideas and constructive discussions on the subject of this paper. CT was supported by the DFG priority program SPP 2265 \textit{Random Geometric Systems}.

\end{document}